\documentclass[11pt]{article}
\usepackage[mathvars]{brsheader} 
\abbreviateBBMathFont

\title{A canonical Euclidean Gallai theorem}
\author{Benedict \textsc{Randall Shaw}\footnotemark[1]}
\date{November 2026}

\begin{document}

\maketitle

\renewcommand{\thefootnote}{\fnsymbol{footnote}}
\footnotetext[1]{\href{mailto:bwr26@cam.ac.uk}{bwr26@cam.ac.uk}, Department of Pure Mathematics and Mathematical Statistics (DPMMS), University of Cambridge, Wilberforce Road, Cambridge, CB3 0WA, United Kingdom}

\begin{abstract}
In this note, we prove a canonical Euclidean Gallai theorem: for any finite set of points \(S\subset \mathbb{R}^d\), any colouring of \(\mathbb{R}^d\) contains a scaled (and possibly rotated) copy of \(S\) that is either monochromatic or rainbow.
\end{abstract}

\section{Introduction}
We are concerned with colourings of \(\mathbb{R}^d\)---that is, maps from \(\mathbb{R}^d\) to some other set, whose structure we neglect, and which for our purposes may be infinite. Equivalently, then, we may think of such a colouring as an equivalence relations on \(\mathbb{R}^d\), where two points are equivalent if they receive the same colour. We call a colouring using finitely many colours a \emph{finite colouring}, and a colouring using at most \(r\) colours an \(r\)-\emph{colouring}. Under a given colouring, we say a set of points is \emph{monochromatic} if all the points receive the same colour, or \emph{rainbow} if all the points receive different colours.

Euclidean Ramsey theory has historically primarily investigated finite colourings. One of the main areas of research was the question of exactly which sets are \emph{Ramsey}---that is, which sets \(S\) have the property that for any \(r\), we find that for sufficiently large \(d'\) all \(r\)-colourings of \(\mathbb{R}^{d'}\) contain a congruent monochromatic copy of \(S\). This was introduced by Erd\H{os}, Graham, Montgomery, Rothschild, Spencer, and Straus \cite{erdos1973}, saw major contributions including Frankl and R\"odl \cite{frankl1990}, Kříž \cite{kriz1991}, Leader, Russell, and Walters \cite{leader2010transitivesetseuclideanramsey}, and was ultimately closed by OpenAI \cite{openai2026}. If one instead considers homothetic copies---those that are scaled and translated, but not rotated---we have a much cleaner picture of the world in Gallai's theorem (first mentioned in \cite{rado1945}).

\begin{theorem}[Gallai]
    For any finite set of points \(S\subset \mathbb{R}^d\), any finite colouring of \(\mathbb{R}^d\) contains a homothetic copy of \(S\) that is monochromatic.
\end{theorem}

Of course, if we remove the restriction that the set of colours be finite, we can no longer hope to find a monochromatic copy of any configuration, since we could very well just colour every point differently. Instead, we seek a copy which is either monochromatic or rainbow. There has been recent interest in the question of which sets are \emph{canonically Ramsey}---that is, for which sets \(S\) some \(d'\) has the property that any colouring of \(\mathbb{R}^{d'}\) contains a congruent copy of \(S\) that is either monochromatic or rainbow. This was introduced by Mao, Ozeki, and Wang \cite{mao2022euclideangallairamseytheory}, saw incremental progress by Cheng and Xu \cite{cheng2023euclideangallairamseyvariousconfigurations}, Geh\'er, Sagdeev, and T\'oth \cite{geher2026}, and Fang, Ge, Shu, Xu, Xu, and Yang \cite{fang2025canonicalramseytrianglesrectangles}. This year saw wider classes of canonically Ramsey sets established by Ge Shu, Xu, and Yu \cite{ge2026simplicesexhibitcanonicalramsey} and the author \cite{randallshaw2026cuboidscanonicallyramsey,randallshaw2026productssimplicescanonicallyramsey,randallshaw2026regularpentagoncanonicallyramsey}, until Sagdeev \cite{sagdeev2026ramseysetscanonicallyramsey} gave a proof using ChatGPT that all Ramsey sets are canonically Ramsey, and gave an example of a canonically Ramsey set that is not itself Ramsey.

In the setting of potentially infinite colourings, we can no longer ask for a homothetic copy to be either monochromatic or rainbow---for example, colouring \(\mathbb{R}^2\) by \(x\)-coordinate induces a colouring in which every homothetic copy of \(\{0,1\}^2\) receives exactly two colours. We therefore must also allow our scaled copies to be rotated, and so our main result is the following theorem:

\begin{theorem}\label{thm:main}
    For any finite set of points \(S\subset \mathbb{R}^d\), any colouring of \(\mathbb{R}^d\) contains a scaled (and possibly rotated) copy of \(S\) that is either monochromatic or rainbow.
\end{theorem}

\section{Proof}

We will use the following canonical version of Gallai's theorem on \(\mathbb{N}^k\), due to Deuber, Graham, Pr\"omel, and Voigt \cite{Deuber1983ACP}, with an elementary proof later given by Pr\"omel and R\"odl \cite{promel1986}:

\begin{theorem}[\cite{Deuber1983ACP}]\label{thm:Ncanonicalgallai}
    For any \(m,k\), there is some \(n\) such that, within any colouring of \([n]^k\), there is some linear subspace \(U\) of \(\mathbb{Q}^k\), and some set of points of \([n]^k\) homothetic to \([m]^k\), within which any two points \(x,y\) have the same colour if and only if \(x-y\in U\).
\end{theorem}

Note several crucial differences between this theorem and our result: first, our result also applies to sets of points which are not similar to a subset of \(\mathbb{N}^d\); and secondly, we find a copy with one of two possible colourings---monochromatic or rainbow---rather than the wider class of colourings permitted by Theorem \ref{thm:Ncanonicalgallai}. Of course, this comes at the cost of also including copies which have been rotated as well as scaled.

\begin{proof}[Proof of Theorem \ref{thm:main}]
Say the set of points we are concerned with has \(t\) points \(S=\{a_1,\dots,a_t\}\). We will choose \(k=t+t{t\choose 2}\left({t\choose 2}-1\right)\), and identify \([k]\) with the set \([t]\) together with sets of the form \[(p,q,i): p,q\in [t]^{(2)}, p\neq q, i\in [t].\]
In \(\mathbb{Q}^k\), therefore, we may refer to basis vectors \(\mathbf{e}_i\) and \(\mathbf{e}_{p,q,i}\). We now build a linear map \(\varphi: \mathbb{N}^k\to \mathbb{R}^d\), by choosing its images on these basis vectors such that:---
\begin{itemize}
    \item each \(\mathbf{e}_i\) is mapped to \(a_i\),
    \item each \(\mathbf{e}_{p,q,i}\) is mapped to some \(a_{p,q,i}\), such that for each \((p,q)\), some similarity maps each \(a_i\) to the corresponding \(a_{p,q,i}\), and
    \item for each \(p,q\),
    \[\{a_i: i\in p\} = \{a_{p,q,i}: i \in q\}.\]
    That is, the copy of \(S\) comprising the points \(a_{p,q,i}:i\in [t]\) sits on top of \(S\) itself in such a way that the points indexed by \(q\) are exactly those of \(S\) itself which are indexed by \(p\).
\end{itemize}
Note that this last requirement will require our copies of the form \(\left\{a_{p,q,i}:i\in [t]\right\}\) to be merely similar to \(S\), rather than homothetic. We may extend our notation to consider \(a_{p,p,i}=a_i\) itself, noting that this trivially coincides with the points of \(S\) in the correct way, and likewise \(\mathbf{e}_{p,p,i}=\mathbf{e}_i\).

We extend \(\varphi\) to the whole of \(\mathbb{N}^k\), and notice that we may now consider the colouring of \(\mathbb{N}^k\) pulled back from our colouring of \(\mathbb{R}^d\), whereby two points of \(\mathbb{N}^k\) receive the same colour if and only if their images under \(\varphi\) do. Notice that by the definition of the \(a_{p,q,i}\), if for some \(p=\{i,j\}\) \(\mathbf{e}_i\) and \(\mathbf{e}_j\) receive different colours, then so too for any \(q={k,l}\) do \(\mathbf{e}_{p,q,k}\) and \(\mathbf{e}_{p,q,l}\).

We now choose \(m\geq {t \choose 2}\left({t\choose 2}+1\right)\), and apply Theorem \ref{thm:Ncanonicalgallai} to find a linear subspace \(U\) of \(\mathbb{Q}^k\), and a homothetic copy of \(\{0,\dots,m\}^k\) in which any two points receive the same colour if and only if they differ by an element of \(U\). Restricting \(\varphi\) to this copy of \(\{0,\dots,m\}^k\), we have a map \(\varphi':\{0,\dots,m\}^k\to \mathbb{R}^d\) in, just as before:---
\begin{itemize}
    \item each \(\mathbf{e}_i\) is mapped to some \(a'_i\), such that some homothety maps each \(a_i\) to \(a'_i\)
    \item each \(\mathbf{e}_{p,q,i}\) is mapped to some \(a'_{p,q,i}\), such that for each \((p,q)\), some similarity maps each \(a'_i\) to the corresponding \(a'_{p,q,i}\), and
    \item for each \(p,q\),
    \[\{a'_i: i\in p\} = \{a'_{p,q,i}: i \in q\}.\]
\end{itemize}
In particular, for any \(x\in \{0,\dots,m\}^k\), and choice of \(b\) and of \(c_{p,q}\) for each \(p,q \in [t]^{(2)}, p \neq q\), the set of points
\[\left\{\varphi'\left(b\mathbf{e}_i +  \sum_{p,q}c_{p,q}\mathbf{e}_{p,q,i}\right):i\in t\right\}\]
is similar to \(S\), unless they are in fact all the same point. So it will suffice to find some \(b, c_{p,q}\) on which
\(\left\{b\mathbf{e}_i +  \sum_{p,q}c_{p,q}\mathbf{e}_{p,q,i}:i\in [t]\right\}\) is monochromatic or rainbow in \(\{0,\dots,m\}^k\)---noting that if we find these points are monochromatic, we will also need to check that they are distinct.

We are certainly done if all \(\mathbf{e}_i-\mathbf{e}_j\) are in \(U\), as then \(\mathbf{e}_1,\dots,\mathbf{e}_t\) is monochromatic, and these certainly correspond to distinct points in \(\mathbb{R}^d\). So assume otherwise: that some \(p=\{i,j\}\in [t]^{(2)}\) has \(\mathbf{e}_i-\mathbf{e}_j\notin U\). Thus, recalling our previous identification of \(\mathbf{e}_{p,p,i}\) with \(\mathbf{e}_i\), we notice that for each \(q={k,l}\), \(\mathbf{e}_{p,q,k}-\mathbf{e}_{p,q,l}\notin U\). 

We now set \(c_{p',q}=0\) for each \(p\neq p'\), and choose each \(c_{p,q}\), as well as \(b\), which we consider to be \(c_{p,p}\), uniformly randomly from \({t\choose 2} + 1\). Notice that as \(m\geq {t \choose 2}\left({t\choose 2}+1\right)\), we certainly have
\[\left\{\sum_{q\in[t]^{(2)}}c_{p,q}\mathbf{e}_{p,q,i}:i\in [t]\right\}\subset \{0,\dots,m\}^k\]
as required. Write \(x_i=\sum_{q\in[t]^{(2)}}c_{p,q}\mathbf{e}_{p,q,i}\). It will now suffice to show that for each pair \(i,j\), the probability that \(x_i\) and \(x_j\) receive the same colour is at most \(\left({t\choose 2} + 1\right)^{-1}\)---as then there will be a probability of at least \(\left({t\choose 2} + 1\right)^{-1}\) that this does not happen for any \(i,j\), and so the set of \(x_i\) is rainbow, as desired.

Since the \(c_{p,q}\) for each \(q\) are chosen independently at random, we may assume that for all \(q\neq \{i,j\}\), we have already fixed some \(c_{p,q}\). So it will be enough to show that, given all other \(c_{p,q}\), at most one of the \({t\choose 2} + 1\) possible choices of \(c_{p,\{i,j\}}\) causes \(x_i\) and \(x_j\) to have the same colour. Write \(x'_\ell=\sum_{q\in[t]^{(2)}\setminus \{\{i,j\}\}}c_{p,q}\mathbf{e}_{p,q,\ell}\), so that \(x_\ell=x'_\ell + c_{p,\{i,j\}}\mathbf{e}_{p,\{i,j\},\ell}\). Then \(x_i\) and \(x_j\) will have the same colour only if \[(x'_i - x'_j) + c_{p,\{i,j\}}\left(\mathbf{e}_{p,\{i,j\},i}-\mathbf{e}_{p,\{i,j\},j}\right) \in U.\]
But certainly this cannot happen for two different values of \(c_{p,\{i,j\}}\), as then this would contradict the fact that \(\mathbf{e}_{p,\{i,j\},i}-\mathbf{e}_{p,\{i,j\},j}\notin U\). Hence we have proved the desired result---and so some choice of entire set of \(c_{p,q}\) gives \(x_i\) rainbow, as desired.
\end{proof}

\section{Acknowledgement}
The author is funded by an Internal Graduate Studentship of Trinity College, Cambridge.

	\bibliographystyle{plain}
\bibliography{main}

\end{document}